\documentclass[sn-mathphys-num]{sn-jnl}

\usepackage{graphicx}%
\usepackage{multirow}%
\usepackage{amsmath,amssymb,amsfonts}%
\usepackage{amsthm}%
\usepackage{mathrsfs}%
\usepackage[title]{appendix}%
\usepackage{xcolor}%
\usepackage{textcomp}%
\usepackage{manyfoot}%
\usepackage{booktabs}%
\usepackage{algorithm}%
\usepackage{algorithmicx}%
\usepackage{algpseudocode}%
\usepackage{listings}%

\theoremstyle{thmstyleone}%
\theoremstyle{thmstyletwo}%

\theoremstyle{thmstylethree}%

\usepackage{fullpage} 
\usepackage{amsmath}
\usepackage{amssymb}
\usepackage{amsthm}
\usepackage{graphicx}
\newcommand{\concatA}{%
  \mathbin{\raisebox{1ex}{\scalebox{.7}{$\frown$}}}%
}
\newcommand{\concatB}{%
  \mathbin{\rotatebox[origin=c]{90}{\scalebox{.7}{(\kern1ex)}}}
}

\newtheorem{thm}{Theorem}[section]
\newtheorem{lem}[thm]{Lemma}

\newtheorem{qu}{Question}

\theoremstyle{definition}

\newtheorem{defn}{Definition}

\newtheorem{rmk}{Remark}

\newcommand{\setof}[2]{\lbrace #1\, | \, #2 \rbrace}

\begin{document}

\title[Article Title]{A closed subset of Baire space not Medvedev equivalent to any closed set of Cantor space }


\author{\fnm{Joshua A.} \sur{Cole}}\email{jcole@benedictine.edu}



\affil{\orgdiv{Department of Mathematics and Computer Science}, \orgname{Benedictine College}, \orgaddress{\street{1020 North Second Street}, \city{Atchison}, \postcode{66002}, \state{KS}, \country{USA}}}




\abstract{For \emph{mass problems} $P,Q\subseteq {\mathbb{N}^\mathbb{N}}$ (\emph{Baire space}), $P$ is \emph{Medvedev reducible} to $Q$ ($P\leq_sQ$) if for some Turing funcional $\Phi$, $\Phi(Q)\subseteq P$, and \emph{Medvedev equivalent} to $Q$ if also $Q\leq_sP$. Shafer asked if every closed problem $P$ is Medvedev equivalent to a closed problem $Q$ with $Q\subseteq 2^\mathbb{N}$ (\emph{Cantor space}). We show that this is not the case.}




\maketitle

\section{Introduction}\label{sec1}

\begin{defn} A \emph{mass problem}  is a subset $P$ of $\mathbb{N}^{\mathbb{N}}$, the set of all functions from the set $\mathbb{N}$ of natural numbers into itself. An element $f \in P$ is called a \emph{solution} to $P$. There are two natural ways to compare mass problems. For any mass problems $P$ and $Q$,
\begin{enumerate}
\item $P\leq_wQ$ ($P$ is \emph{weakly} or \emph{Muchnik reducible} to $Q$) iff for every $g\in Q$, there is some $f\in P$ such that $f\leq_Tg$ ($f$ is \emph{Turing reducible to} or \emph{computable in} $g$);
\item $P\leq_sQ$ ($P$ is \emph{strongly} or \emph{Medvedev reducible} to $Q$) iff there is a Turing functional $\Phi$ such that for every $g\in Q$, $\Phi(g)\in P$; we also write $\Phi(Q)\subseteq P$.
\item for $\bullet$ either $w$ or $s$, $P\equiv_\bullet Q$ ($P$ and $Q$ are \emph{Muchnik (Medvedev) equivalent}) iff both $P\leq_\bullet Q$ and $Q\leq_\bullet P$.
\end{enumerate}
\end{defn}
Since both of these relations are \emph{preorderings} (transitive and reflexive) but in general two elements need not be comparable, they are often converted to \emph{partial orderings} of the equivalence classes of mass problems defined by mutual reducibility $\equiv_\bullet$. The resulting partially ordered structures are known as the \emph{Muchnik (Medvedev) degrees}.

Each of the sets $\mathbb{N}^{\mathbb{N}}$ and $2^{\mathbb{N}}$ is naturally considered as a topological space with a topology generated by the sets $[\sigma] := \setof{f}{\sigma\subseteq f}$, where $\sigma\in \mathbb{N}^{<\mathbb{N}}$ is a finite sequence of natural numbers and $\sigma\subseteq f$ means that $\sigma$ is a finite initial segment of $f$ viewed as an infinite sequence. $\mathbb{N}^{\mathbb{N}}$ is commonly known as \emph{Baire space} while $2^{\mathbb{N}}$ is known as \emph{Cantor space}. One of the important differences between the two, which is particularly relevant for us here, is that Cantor space is a compact space while Baire space is not.

All open sets in this topology have computable elements and so are trivial in the sense that they lie below all other sets in both of $\leq_w$ and $\leq_s$ and are thus all mutually equivalent. Hence the simplest topologically characterized sets about which there are interesting questions about either reducibility are the closed sets. Questions of this sort fit into the general category of comparisons between two measures of complexity that is a frequent topic for computability theorists.

Paul Shafer \cite{Sh2} (Question 1.4) has posed questions about how the theory of closed mass problems  in Baire space and preordered by Muchnik or Medvedev reducibility compares to the theory of closed mass problems in Cantor space.  In particular he asked whether every closed subset of Baire space is Medvedev equivalent to some closed subset of Cantor space.  In this paper, we provide a negative answer (Theorem \ref{thm:main}).

Other authors have also studied how topologically-defined classes behave as mass problems. For example Lewis, Shore, and Sorbi \cite{LeShSo} studied dense, closed, and discrete Medvedev degrees; Shafer \& Sorbi \cite{ShSo} studied the relationship between the closed Medvedev degrees and the degrees of enumerability.

Another direction of research that has been productive is the study of effective versions of topological classes, most prominently the $\Pi^0_1$ or \emph{effectively closed} mass problems in Cantor space. These are the main focus of Simpson \cite{Si}, which is also a good general introduction to mass problems. For a wide-ranging survey of results on Muchnik and Medvedev reducibility, see Hinman \cite{Hi1}.  Like Simpson, Hinman focuses more on $\Pi_1^0$ subsets of $2^{\mathbb{N}}$, especially in deciding which proofs to include, but does cover other areas as well. The bibliographies in these two surveys are excellent guides to the literature of mass problems.

\section{Notation  \& Background} 
\label{sec:background}

In the following, $\sigma, \tau \in \mathbb{N}^{<\mathbb{N}}$ -- finite sequences of natural numbers. We write $\sigma \subset \tau$ to mean that $\sigma$ is a proper initial segment of $\tau$, and as above, $\sigma \subset f$ means $\sigma$ is an initial segment of $f$. We use $\langle \rangle$ to denote the empty sequence. $|\sigma|$ is the length of $\sigma$. By $\tau \concatA \sigma$ is meant the concatenation of $\tau$ with $\sigma$. If $m \in \mathbb{N}$, then $\sigma \concatA m$ is the concatenation of $\sigma$ with the sequence whose only entry is $m$. A \emph{tree} is a set finite sequences closed under subsequence. If $T$ is a tree, then $[T]$ is the set of all infinite paths through $T$; for a set $S$ not necessarily closed under subsequence we write $S^*$ for its closure and $[S]$ for $[S^*]$. If $P$ is a $\Pi_1^0$ subset of Baire space, then $T_P$ denotes a computable tree such that $[T_P] = P$. 

If $\Phi$ is a Turing functional, then $\Phi(\sigma) = \tau$ just in case $\tau$ is maximal such that for all $n \le |\tau|, \Phi_{|\sigma|}^{\sigma}(n) = \tau(n)$, where the subscript $|\sigma|$ means that in some canonical formalism the computation takes place in at most $|\sigma|$ steps. Of course, it is possible that $\Phi(\sigma) = \langle \rangle$. 

\begin{defn} 
\label{dfn:ordering}
For $\sigma \in \mathbb{N}^{<\mathbb{N}}$, we define $\displaystyle \#(\sigma) = \prod_{i < |\sigma|} p_i^{\sigma(i)}$, where $p_0, p_1, p_2, \cdots$ is the enumeration of the prime numbers in increasing order. Note that this provides an ordering of type $\omega$ on $\mathbb{N}^{<\mathbb{N}}$.

\end{defn}

Other notation is standard in computability theory, as in, for example, Soare \cite{So}. 

Finally in this section:

\begin{lem}[Simpson]\label{lem:CS}(See Cole/Simpson \cite{CoSi}, proof of Sublemma 5.6.)
 For any $\Pi^0_2$  set $P\subseteq\mathbb{N}^{\mathbb{N}}$, there exists a  $\Pi^0_1$ set $\hat{P}\subseteq\mathbb{N}^{\mathbb{N}}$ and a computable homeomorphism $\Phi$ mapping $P$ onto $\hat{P}$. In particular $P\equiv_s\hat{P}$.
\end{lem}
\begin{proof}
By the definition of $\Pi^0_2$ there exists a computable relation $R$ such that
$$P=\setof{f}{\forall m\exists nR(f,m,n)}.$$
For $f\in P$, set 
\begin{align*}
g_f(m) :&= \text{ least } n \text{ such that } R(f,m,n),\\
\Phi(f):&=f\oplus g_f \text{ and}\\
\hat{P}:&=\setof{f\oplus g_f}{f\in P}\\
&= \setof{f\oplus g}{\forall m [R(f,m,g(m)) \land (\forall n<g(m))\lnot R(f,m,n)]}.
\end{align*}
Then the conclusion of the Lemma is clear.
\end{proof}

\section{A closed set in Baire Space not Medvedev equivalent to any closed set in Cantor Space}
\label{sec:closed}

We will deduce our result from the following theorem.

\begin{thm}
\label{thm:main}
 There is a nonempty closed (in fact $\Pi_1^0$) subset of Baire space $Q$ that is not compact and such that if $f,g \in Q$ are distinct, then $f \nleq_T g$; that is, any two distinct elements of $Q$ are Turing incomparable. 
 \end{thm}

A version of this theorem for Cantor space was a classic result of Jockusch and Soare (\cite{JoSo}, Theorem 4.7).

First we show how this theorem answers Shafer's question in the negative. Then we'll prove the theorem by a priority argument inspired by the original argument of Jockusch and Soare, but with significant alterations. (Their finite-injury argument cannot be transferred in a straightforward way, because it makes essential use of the fact that for any $n \in \mathbb{N}$, there are only finitely many binary strings  (elements of $2^{<\mathbb{N}}$) of length $n$, which is false for  $\mathbb{N}^{<\mathbb{N}}$.)

\begin{thm}[Main Result] There is a closed (in fact $\Pi_1^0$) subset of Baire space not Medvedev (strongly) equivalent to any closed subset of Cantor space. 
\end{thm}

\begin{proof}
Let $Q$ be a subset of Baire space as given by Theorem \ref{thm:main} and suppose towards a contradiction that $Q$ is Medvedev equivalent to a closed subset $P$ of Cantor space. Then there exit Turing functionals $\Phi$ and $\Psi$ such that $\Phi(Q)\subseteq P$ and $\Psi(P)\subseteq Q$. From basic topology $P$ as  a closed subset of a compact space is itself compact and thus that $\Psi(P)$ as a continuous image of a compact set is also compact. Since $Q$ is not compact, $\Psi(P)\neq Q$, so there exists some $f\in Q/\Psi(P)$. Now $\Phi(f)\in P$, so $g:=\Psi(\Phi(f))\in \Psi(P)$.  Also $g\in Q$. Since both $f$ and $g$ belong to $Q$ and clearly $g\leq_Tf$, by the defining property of $Q$, $f=g$. But then also $f\in\Psi(P)$, contrary to the choice of $f$. 
\end{proof}

We turn now to:

\begin{proof}(of Theorem \ref{thm:main})
As mentioned above, this construction will be in the style of Jockusch/Soare \cite{JoSo}; however its form and notation are also influenced by Binns/Simpson \cite{BiSi} and by 
\newline
Binns/Shore/Simpson \cite{BiShSi}.

Our priority construction builds a sequence of maps from finite strings to finite strings, from which will be derived the desired non-compact $\Pi_1^0$ subset $Q$ of Baire space whose members are pairwise Turing incomparable. 

\begin{defn}[Simpson] A \emph{treemap} is a function $F: \mathbb{N}^{<\mathbb{N}} \rightarrow \mathbb{N}^{< \mathbb{N}}$ such that   \[  F(\sigma^\frown i) \supseteq F(\sigma)^\frown i, \] for all $i \in \mathbb{N}$ and $\sigma \in \mathbb{N}^{<\mathbb{N}}$.
\end{defn}

\begin{lem} Suppose $F: \mathbb{N}^{<\mathbb{N}} \rightarrow \mathbb{N}^{< \mathbb{N}}$ is a treemap. 
\label{lem:treemap1}
\begin{enumerate}

\item $F(\sigma^\frown i)$ and $F(\sigma^\frown j)$ are incompatible extensions of $F(\sigma)$ for distinct $i,j \in \mathbb{N}$.
\item Each treemap $F$ determines in a uniformly effective way a tree \[F(T) = \lbrace \tau \in \mathbb{N}^{\mathbb{N}} \, | \, \exists \sigma, F(\sigma) \supseteq \tau \rbrace. \]  Moreover, for each $\tau \in F(T)$ the least such $\sigma$ as in the definition of $F(T)$ will be a subsequence of $\tau$ (Cole/Simpson \cite[Remark 5.3]{CoSi}). Therefore the quantifier in the definition of $F(T)$ can be bounded.
\item From part two of this Lemma, it follows that if a treemap $F$ is computable, so is the tree $F(T)$, via a computation uniform in $F$.

\end{enumerate}

\end{lem}
\begin{proof} Immediate from the definition of a treemap.
\end{proof}

\begin{lem} 
\label{lem:string}
Suppose $F$ is a treemap  and $F(T)$ is the tree defined in part (2) of Lemma ~\ref{lem:treemap1} above. Then for each $f \in [F(T)]$, for every $n \in \mathbb{N}$, there is a finite string $\sigma \in \mathbb{N}^{<\mathbb{N}}$ such that $f \upharpoonright n \subseteq F(\sigma) \subset f$. 

\end{lem}

\begin{proof} Immediate from the definitions. 
\end{proof}

\begin{defn} 
\label{defn:nestedsequence}
A \emph{nested sequence of tree maps} is a sequence $\langle F_s \rangle_{s \in \mathbb{N}}$ of treemaps $F_s: \mathbb{N}^{<\mathbb{N}} \rightarrow \mathbb{N}^{< \mathbb{N}}$ that satisfy the following three properties. 
\begin{enumerate}
\item $F_0(\sigma) = \sigma$ for all $\sigma \in \mathbb{N}^{<\mathbb{N}}$.

\item range($F_{s+1}) \subseteq$ range($F_s$), for all $s \in \mathbb{N}$.
 \item $F(\sigma) = \lim_s F_s(\sigma)$ exists for all $\sigma \in \mathbb{N}^{<\mathbb{N}}$. 
\end{enumerate}
\end{defn}

\begin{rmk}  As Jockusch and Soare note, the use in priority constructions of nested sequences of tree maps goes back to Shoenfield \cite{Sh}.  \end{rmk}

\medskip

Our construction builds a $\mathbf{0'}$-computable nested sequence of \textbf{computable} treemaps, $\langle F_s \rangle_{s \in \mathbb{N}}$. We repeat that each individual treemap is computable, but we make use of a $0'$ oracle to compute an index for each $F_s$ as a partial-computable function. 

For each treemap $F_s$, let $T_s$ be the tree $F_s(\mathbb{N}^{\mathbb{N}})$ defined as indicated in Lemma ~\ref{lem:treemap1}. Since the sequence $\langle F_s \rangle_{s \in \mathbb{N}}$ is $0'$-computable, by Lemma \ref{lem:treemap1} part (3)  the sequence $\langle T_s \rangle_{s \in \mathbb{N}}$ is a $0'$-computable sequence of computable trees.

After the construction, we'll define $\hat{Q} = \bigcap_{s \in \mathbb{N}} [T_s]$. The definition of a treemap and Property (3) of the definition of a nested sequence of treemaps guarantee that  $\hat{Q}$ is nonempty. As each $[T_s]$ is closed, $\hat{Q}$ is closed, which is all we need for our main theorem. However, we go even further in our analysis: \[ f \in \hat{Q} \iff \forall s (f \in [T_s]) \iff \forall s \forall n 
(f \upharpoonright n \in T_s). \]

Since the sequence $\langle T_s \rangle_{s \in \mathbb{N}}$ is a $0'$-computable sequence of computable trees, $\hat{Q}$ is a $\Pi_1^{0, 0'}$ class. By Post's Theorem, it follows that $\hat{Q}$ is a $\Pi_2^0$ class. By Lemma \ref{lem:CS} there is a $\Pi_1^0$ subset of Baire space $Q$ homeomorphic to $\hat{Q}$ via a  homeomorphism that preserves Turing degree. 

By properties (1) and (3) of the definition of a nested sequence of treemaps (Definition ~\ref{defn:nestedsequence}), $\hat{Q}$ is homeomorphic to Baire space, and so $\hat{Q}$ is non-compact. Thus the aforementioned $\Pi_1^0$ class $Q$ homemorphic to $\hat{Q}$ will be non-compact, one of the properties needed to satsify Theorem ~\ref{thm:main}. The main work of the construction will be to ensure that distinct elements of $\hat{Q}$, and hence of $Q$, are pairwise incomparable in Turing degree. 

\medskip
 For each $\sigma \in \mathbb{N}^{<\mathbb{N}}$ and each $e \leq |\sigma|$ our construction has the requirement:

\[ R_{\sigma,e} \equiv \text{ for all } f \text{ extending } F(\sigma), \lbrace e \rbrace^f \supseteq F(\sigma),    \text{ or } \lbrace e \rbrace^f \notin \hat{Q}. \]

\medskip
$R_{\sigma,e}$ is \textbf{satisfied at stage} $\mathbf{s}$ if $\lbrace e \rbrace^{F_s(\sigma)}$ is comparable with  $F_s(\sigma)$ or $\lbrace e \rbrace^{F_s(\sigma)} \notin T_s$. 

\medskip
We say that $R_{\sigma,e}$ \textbf{requires attention} at stage $s$ if $\lbrace e \rbrace^{F_s(\sigma)} 
= \nu \in T_s$ and the strings $\nu$ and $F_s(\sigma)$ are incomparable.

\medskip
Recall from Section ~\ref{sec:background}  that $\#$ is a computable function on $\mathbb{N}^{<\mathbb{N}}$ that obeys \[ \sigma \subset \tau \rightarrow \# (\sigma) < \# (\tau). \] If $\# (\sigma) < \#(\tau)$, then any requirement involving $\sigma$ has higher priority than any requirement involving $\tau$. If $i < e$ then $R_{\sigma, i}$ has higher priority than $R_{\sigma, e}$.

\medskip
\textbf{At stage} $\mathbf{0}$ we set $F_0(\sigma) = \sigma$ for all $\sigma \in \mathbb{N}^{<\mathbb{N}}$. Consequently, $T_0 = \mathbb{N}^{<\mathbb{N}}$. We take no action at stage $0$, and the protected list is empty.

\textbf{At stage} $\mathbf{s +1 }$, given $F_s$ and the set of strings on the \emph{protected list}, we consider the highest priority requirement that requires attention, if any. (If none, proceed to the next stage without adding to the protected list and set $F_{s+1}(\sigma) = F_s(\sigma)$ for all $\sigma \in \mathbb{N}^{<\mathbb{N}}$.)  

If $R_{e, \sigma}$ is the highest priority requirement requiring attention, then $\lbrace e \rbrace^{F_s(\sigma)} = \nu \in T_s$, and $F_s(\sigma)$ is incomparable with $\nu$.

\medskip
\begin{defn} 
\label{def:lambda}
At stage $s+1$, if we are considering the requirement $R_{\sigma, e}$ that requires attention with $\nu = \lbrace e \rbrace^{F_s(\sigma)}  \in T_s$, let $\mathbf{\lambda}$ be the longest string such that $F_s(\lambda) \subset \nu$. We call $\lambda$ the \emph{target string} for $R_{\sigma, e}$ at stage $s+1$. 
\end{defn}

Why does such a string $\lambda$ exist when $R_{\sigma, e}$ requires attention? If $\lambda$ does not exist, then there is absolutely nothing in the image of $F_s$ that properly precedes $\nu$. Of course, for $\nu$ to be in $T_s$, $\nu \subseteq F_s(\rho)$ for some string $\rho$. Since $F_s$ is a treemap, we also have $F_s(\langle \rangle) \subseteq F_s(\rho)$. Summing up, for $\lambda$ to not exist, we must have $\nu \subseteq F_s(\langle \rangle)$. But, also because $F_s$ is a treemap and $\langle \rangle \subseteq \sigma$,  it follows that $\nu  \subseteq F_s(\langle \rangle)  \subseteq F_s(\sigma)$, which violates our definition of $R_{\sigma, e}$ requiring attention: $\nu$ and $F_s(\sigma)$ are supposed to be incomparable.

\medskip
We call $\lambda$ the target string because we wish to change its image in order to ensure $\nu \notin T_t$ for all $t > s$.  But we won't be able to act immediately if $\lambda \subset \sigma$ or is on the protected list. A string gets on the \emph{protected list} if we change its image under $F$ at some stage because  it's the target string. The point is, we will allow no string to serve successfully as a target string more than once. 

How might the target string $\lambda$ (with the property $F_s(\lambda) \subset \nu)$ end up as an initial segment of $\sigma$, even if $F_s(\sigma)$ and $\lbrace e \rbrace^{F_s(\sigma)} = \nu$ are incomparable? Here's an example: $\sigma = \lambda^\frown 0$, $\nu = \lambda^\frown1$, and $F_s$ is the identity map. Note that it is impossible that $\lambda \supseteq \sigma$. For that would imply $F_s(\sigma) \subseteq F_s(\lambda) \subset \nu$, but $F_s(\sigma)$ and $\nu$ being incomparable is part of the definition of a requiring attention.
\medskip

We proceed by cases and Case 1 should be thought of as the paradigmatic case. Our actions for other cases are all directed at landing us in the situation of Case 1, if possible. The action taken for Case 1 is just what happens in the above-mentioned constructions of Jockusch/Soare and Binns/Simpson. 

\textbf{Case 1:  The target string $\lambda \not \subset \sigma$ and $\lambda$ is not on the protected list.} (Since as mentioned above $\lambda \supseteq \sigma$ is impossible, when we are in Case 1, $\lambda$ and $\sigma$ are incomparable.) 

\medskip
Pick an $m \in \mathbb{N}$ such that $F_s(\lambda^\frown m)$ is incomparable with $\nu$. Our idea will be to define $F_{s+1}(\lambda) = F_s(\lambda^{\frown} m)$, with the effect that $\nu \notin T_{s+1}$. (It will have been `passed over.') The images of $F$ on successors of $\lambda$ will have to be pushed along accordingly. To be precise:

\[
  F_{s+1}(\gamma) =
  \begin{cases}
                            F_s(\lambda^{\frown}m^{\frown} \gamma')        & \text{ if } \gamma = \lambda^{\frown} \gamma' \text{ for some string } \gamma'. \text{ (Note: } \gamma' \text{ may be } \langle \rangle.) \\
                                   
    F_s(\gamma) & \text{ otherwise. }
  \end{cases}
\]

\begin{rmk}

\noindent
\begin{enumerate}
\item Because $\sigma$ (the string whose requirement $R_{\sigma, e}$ we are acting for) and $\lambda$ are incomparable, under this definition $F_{s+1}(\sigma) = F(\sigma)$. 
\item Since the treemaps will all be nested, this definition of $F_{s+1}$ ensures that \[ \lbrace e \rbrace^{F_{s+1}(\sigma)} \notin T_{t} \text{ for all } t > s.  \]
\item If $F_t(\sigma) \supseteq F_{s+1}(\sigma)$ for all $t \geq s+1$, the first two remarks here ensure that $R_{\sigma, e}$ will be satisfied. 
\item  $F_{s+1}(\gamma^\frown i) \supseteq F_{s+1}(\gamma)^\frown i$, for all $i \in \mathbb{N}$ and $\gamma \in \mathbb{N}^{\mathbb{N}}$, assuming this was true at stage $s$. (The proof is by induction on the length of $\gamma$.)

\end{enumerate}

\end{rmk}

Our final action in Case 1 is to add $\lambda$ to the \emph{protected list.} The result is that in the future $\lambda$ will never again be used as a target string. 

Then stage $s+1$ ends and we proceed to the next stage.

\medskip
\textbf{Case 2: The target string $\lambda$ for the highest priority requirement $R_{e, \sigma}$ requiring attention is on the protected list or is comparable with $\sigma$.}
\medskip

\emph{Idea behind our treatment of Case 2:} In this case we look (using $0'$) for an extension of $F_s(\sigma)$ in $T_s$ that the Turing functional with index $e$ maps to a proper extension $\nu'$ of $\nu$ in $T_s$. Then we ask if the target string $\lambda'$ associated with $\nu'$  isn't on the protected list and isn't comparable with $\sigma$. Perhaps we find the $\nu'$, but don't immediately have $\lambda'$ meeting both conditions. If this happens, we note our progress and repeat, looking for another extension in the length of agreement. 

We continue until we either fail to find such extensions $\nu'$, or the associated target string $\lambda'$  is not on the protected list and is not comparable with $\sigma$. 

If longer and longer strings $\nu'$ continue to appear, then eventually we must come to a satisfactory $\lambda'$. For long enough $\nu'$ will eventually outgrow deficient $\lambda'$. This is because there are only finitely many initial segments of $\sigma$; and there are only finitely many strings on the protected list at any particular stage.

Suppose that at some point we haven't yet found a suitable target string $\lambda'$ not on the protected list and not comparable with $\sigma$, and we can't push the length of agreement any further (i.e. the `next' $\nu'$ doesn't exist). Then in $F_{s+1}$ we codify our progress in extending the length of agreement as far as we can. If in the future the progress we've codified is never injured, it will be certain we won't see the sort of Turing reduction that would prevent the satisfaction of  $R_{e, \sigma}$. 

To summarize: eventually either the extensions in the length of agreement (i.e. the new $\nu'$) will come to an end, or we will find a suitable target string $\lambda'$.

\medskip
\emph{Case 2 Details:}
Here are the detailed instructions about  about how to proceed if  $\lbrace e \rbrace^{F_s(\sigma)} = \nu \in T_s$, but the target string  $\lambda$ (as in Definition \ref{def:lambda}) is on the protected list or is comparable to $\sigma$. We proceed in substages.

\medskip
At substage $r=0$ we set $\tau_0 = \sigma, \nu_0 = \nu$ and $\lambda_0 = \lambda$.

At substage $r +1$ we ask $0'$:

 \[\text{ (*) Is there is a  }\tau' \supset \tau_r \text{ such that } \lbrace e \rbrace^{F_s(\tau')} = \nu' \in T_s \text{ and } \nu' \supset \nu_r \text{? } \]
 
 \medskip
  \textbf{Case 2A: the answer to (*) is in the \emph{negative}.} We are very happy, and we proceed to define:
\[
  F_{s+1}(\gamma) =
  \begin{cases}
                            F_s(\tau_{r}^{\frown} \gamma')        & \text{ if } \gamma = \sigma^{\frown} \gamma' \text{ for some } \gamma'. \text{ (Note: } \gamma' \text{ may be } \langle \rangle.) \\
                                   
    F_s(\gamma) & \text{ otherwise. }
  \end{cases}
\]

This ends both substage $r+1$ and stage $s+1$ and we proceed to the next stage.

\medskip

\begin{rmk}

\noindent
\begin{enumerate}

\item If $F_{s+1}$ is defined in this way, then no extension of $F_{s+1}(\sigma)$ in $T_s$ can be used as an oracle for the Turing functional indexed by $e$ to compute an element of $[T_s]$ that does not extend $F_{s+1}(\sigma)$. Therefore, provided that $F_t(\sigma) \supseteq F_{s+1}(\sigma)$ for all $t \geq s+1$,  $R_{\sigma, e}$ will be satisfied. 

\item If $f \in [T_{s+1}]$ extends $F_s(\sigma)$, then the first clause of the definition of $F_{s+1}$ ensures that $f$ extends $F_{s+1}(\sigma)$. This is because any string `passed over' by the image of $\sigma$ moving from $F_s(\sigma)$ to 
\newline
$F_{s+1}(\sigma)  \supset F_s(\sigma)$ will not be in $T_{s+1}$. Cf. Jockusch/Soare \cite{JoSo}, last full sentence on page 49.

 \item As in Case 1,  $F_{s+1}(\gamma^\frown i) \supseteq F_{s+1}(\gamma)^\frown i$, for all $i \in \mathbb{N}$ and $\gamma \in \mathbb{N}^{<\mathbb{N}}$, assuming this was true at stage $s$. 

\end{enumerate}
\end{rmk}

\medskip
\textbf{Suppose the answer is yes to the question (*) asked in Case 2 at substage $r+1$ of stage $s+1$.} Then we set $\tau_{r+1}=\tau'$ and $\nu_{r+1} = \nu'$, where $\tau'$ is  least in the ordering given by \# that yields a yes answer, and $\nu'$ is the string that corresponds to $\tau'$ in the question asked by (*).  Then we define $\lambda_{r+1}$ to be the longest $\lambda'$ such that  $F_s(\lambda') \subset \nu'$. 

\medskip
\textbf{Remark:} Note that $\tau_r \supseteq \sigma \rightarrow \tau_{r+1} \supset \sigma$. Since $\tau_0 = \sigma$ it follows that $\tau_{r+1} \supset \sigma$ is always true. 

\medskip
Lastly, we proceed in cases based on whether $\lambda_{r+1}$ is an acceptable target string. 

\medskip
\textbf{Case 2B: The string $\lambda_{r+1}$ is on the protected list or is comparable with $\sigma$.} In this case, simply proceed to the next substage.  

\medskip
\textbf{Claim:} In a  fixed stage $s+1$, Case $2B$ cannot happen infinitely often, since each iteration requires an extension in length of $\nu_r$. 

Specifically:
\begin{itemize}
\item The protected list at stage $s+1$ is finite. 

 \item Considering the problem that $\lambda_{r+1}$ is comparable with $\sigma$, we note that when this happens it must be that $\lambda_{r+1} \subset \sigma$, which will be impossible once $\nu_{r+1}$ is long enough that $\lambda_{r+1}$ ends up longer in length than any initial segment of $\sigma$.  \end{itemize}

The upshot is that if $\nu_{r+1}$ continues to be defined for larger and larger values of $r$, then eventually  $\lambda_{r+1}$ must be incomparable with $\sigma$ and not on the protected list. 

\medskip
\textbf{Case 2C: The string $\lambda_{r+1}$ is not on the protected list and is not comparable with $\sigma$.} Then we define $F_{s+1}$ in a way similar to its manner of definition in Case 1. Pick an $m \in \mathbb{N}$ such that $F_s(\lambda_{r+1}\concatA m)$ is incomparable with $\nu_{r+1}$.

\[
  F_{s+1}(\gamma) =
  \begin{cases}
                    
                                   F_s(\tau_{r+1}\concatA \gamma') & \text { if } \gamma = \sigma^{\frown} \gamma' \text{ for some } \gamma'  \text{ (Note: } \gamma' \text{ may be } \langle \rangle.) \\
                                   
                                           F_s(\lambda_{r+1}\concatA m^{\frown} \gamma')        & \text{ if } \gamma = \lambda_{r+1}\concatA \gamma' \text{ for some } \gamma'  \text{ (Note: } \gamma' \text{ may be } \langle \rangle.) \\
    F_s(\gamma) & \text{ otherwise }
  \end{cases}
\]

\begin{rmk}

\noindent
\begin{enumerate}
\item This is a legitimate definition by cases for the following reason: as $\sigma$ is not comparable with $\lambda_{r+1}$, neither is any extension of one comparable to the other. 
 
 \item This definition ensures that $\lbrace e \rbrace^{F_{s+1}(\sigma)} \notin T_{t}$ for all $t > s$. So long as $F_t(\sigma) \supseteq F_{s+1}(\sigma)$ for all $t \geq s+1$, it follows that $R_{\sigma, e}$ will be satisfied. 

\item As happens if a definition is made in Case 2A, if $f \in [T_{s+1}]$ extends $F_s(\sigma)$, then the first clause of the definition of $F_{s+1}$ ensures that $f$ extends $F_{s+1}(\sigma)$. This is because $\tau_{r+1} \supset \tau_r \supseteq \sigma$, which in turn implies $F_s(\tau_{r+1}) \supset F_s(\sigma)$. Moreover, every extension of $F_s(\sigma)$ not comparable with $F_{s+1}(\sigma)$ will be missing from $T_{s+1}$.

\item As in Case 1 and Case 2A, $F_{s+1}(\gamma^\frown i) \supseteq F_{s+1}(\gamma)^\frown i$, for all $i \in \mathbb{N}$ and $\gamma \in \mathbb{N}^{<\mathbb{N}}$, assuming the stage $s$ version of this fact was true.
\end{enumerate}
\end{rmk}

\medskip
After making this definition of $F_{s+1}$, we place $\lambda_{r+1}$ on the protected list. The result is that in the future $\lambda$ will never again be used as a target string.

This ends both substage $r+1$ and stage $s+1$, and we proceed to the next stage.

\medskip
This completes the construction.

\medskip
\noindent
\textbf{Verification:}

\begin{lem}
\label{lem:ext}
If $F_{s+1}(\gamma) \neq F_s(\gamma)$ for some string $\gamma \in \mathbb{N}^{<\mathbb{N}}$, one of the following holds:
\begin{enumerate}
\item $\gamma$ is a target string for which action is taken at stage $s+1$. 
\item A definition of $F_{s+1}$ is made in Case 2A or Case 2C for the sake of some requirement $R_{\gamma, e}$ with $e \leq |\gamma|$.
\item  An initial segment $\gamma^* \subset \gamma$ is a target string for which action is taken at stage $s+1$.
\item A definition of  $F_{s+1}$ is made in Case 2A or Case 2C  for some requirement $R_{\gamma^*, e}$, where $\gamma^* \subset \gamma$ and $e \leq |\gamma^*|$.
\end{enumerate}

If (1) or (2) is true, then $F_{s+1}(\gamma) \supset F_s(\gamma)$. Thus, if (1) or (2) is true, we don't consider the action taken in step $s+1$ to have injured any requirement $R_{\sigma, e}$. 
\end{lem}
\begin{proof}
By inspection of the possible definitions of $F_{s+1}$. 
\end{proof}

\begin{rmk}
If 3 or 4 is the case, then it is very likely (but not certain) $F_{s+1}(\gamma)$ is incomparable with $F_s(\gamma)$. The exception to this last claim occurs when $F_s(\gamma)$ was already on the new route to be taken by $F_{s+1}(\gamma^*)$. 
\end{rmk}

\begin{defn} A requirement $R_{\gamma, e}$ is \emph{injured} at stage $s+1$ if $F_{s+1}(\gamma) \nsupseteq F_s(\gamma)$.

\end{defn}

\noindent
\begin{rmk} If $R_{\gamma, e}$ is \emph{injured} at stage $s+1$, then (3) or (4) from Lemma \ref{lem:ext} occurred at stage $s+1$. \end{rmk}

\begin{lem}
\label{lem:ind}
For each requirement $R_{\gamma, e}$ action is taken only finitely often. 
\end{lem}
\begin{proof}
We proceed by induction on the priority of $R_{\gamma, e}$. Suppose that action is taken only finitely often for each requirement of higher priority than $R_{\gamma, e}$. 

Then (4) in Lemma \ref{lem:ext} can happen only finitely often by this inductive assumption. Also (3) in Lemma \ref{lem:ext}  can happen only finitely many times, since each string can be used as a target string at most once. Thus $R_{\gamma, e}$ is injured only finitely often. Let $s$ be a stage after which  $R_{\gamma, e}$ is not injured and (by inductive assumption) after which no requirement of higher priority acts.

If $R_{\gamma, e}$ ever requires attention at some stage $t > s$, it will be addressed, since no requirements of higher priority act after stage $s$. After it is addressed, for all $t' \geq  t+1$, $F_{t'}(\gamma) \supseteq F_{t+1}(\gamma)$, since $R_{\gamma, e}$ is never injured after stage $s$.  And then by Remark 2.3, Remark 3.1, and Remark 4.2, $R_{\gamma,e}$ will remain satisfied at each stage $t' \geq t+1$. Hence $R_{\gamma, e}$ will never want attention at any stage $t' > t$, and will therefore never act after stage $t$, and so acts only finitely often.

  \end{proof}
  
\begin{lem}
\label{lem:finite}
For each string $\gamma$, for only finitely many stages $s$ is it true that $F_{s+1}(\gamma) \neq F_s(\gamma)$. 
\end{lem}
\begin{proof}
This follows from Lemma \ref{lem:ind} and its proof. In particular, let $e = |\gamma|.$ Once we've reached a stage after which no action is taken for any requirement of higher priority than $R_{\gamma, e}$ or with some $\gamma^* \subseteq \gamma$ as target string, $R_{\gamma, e}$ itself acts at most once (say at stage $t$) and $R_{\gamma, e}$ is satisfied at each stage after $t$. Then at any later stage, no action of any requirement will change the image of $\gamma$: for all  $t' \geq t+1$, $F_{t'}(\gamma) = F_{t+1}(\gamma)$. 

\end{proof}

\begin{rmk} Lemma \ref{lem:finite} implies $F(\gamma) = \lim_t F_t(\gamma)$ exists for all strings $\gamma$. This fact ensures that $\hat{Q}$ is homeomorphic to Baire space (and so is not compact). \end{rmk}

Recall $\hat{Q} = \lbrace f \in \mathbb{N}^{\mathbb{N}}\, | \,  \forall s \forall n (f \upharpoonright n \in T_s) \rbrace$, where $T_s = \lbrace \tau \, | \, \exists \sigma, F_s(\sigma) \supseteq \tau \rbrace.$

\medskip

\textbf{Claim:} If $f \neq g$ are in $\hat{Q}$, then $f, g$ are Turing incomparable. 
\begin{proof} Suppose by way of contradiction that $\lbrace e \rbrace^f = g$,  $f \neq g$, and $f, g \in \hat{Q}$. Let $n$ be such that 
 $f \upharpoonright n \neq g \upharpoonright n$. Let $\sigma, s$ be such that 
 \begin{itemize}
\item $|\sigma| \geq e$ 
\item $f \upharpoonright n \subseteq F_s(\sigma) \subset f$
\item $\lbrace e \rbrace^{F_s(\sigma)} \supseteq g \upharpoonright n$
\item $R_{\sigma, e}$ is not injured after stage $s$, and 
\item no requirement of higher priority than $R_{\sigma, e}$ acts after stage $s$. 
\end{itemize}

The existence of such $\sigma, s$ is justified by the assumptions on $f$ \and $g$, Lemma ~\ref{lem:string}, Lemma \ref{lem:ind}, and the fact that each string can be the target string at most once.  From these facts it follows that $R_{\sigma, e}$ requires and receives attention at stage $s+1$. Since $\lbrace e \rbrace^f = g$, by Remark 3.1 after the definition made in Case 2A, the action taken in stage $s $ will result from Case 1 or from Case 2C. In either case, the result is that $\lbrace e \rbrace^{F_{s+1}(\sigma)} \notin T_{t}$ for all $t \geq s+1$.  By Remark 4.3 after the definition in Case 2C,  $F_{s+1}(\sigma) \subset f$. (This is trivial in Case 1, since there $F_{s+1}(\sigma)=F_s(\sigma)$, which by the definition of $s$ here is an initial segment of $f$.)  Thus $\lbrace e \rbrace^f \notin \hat{Q}$, contradicting that $\lbrace e \rbrace^f = g \in \hat{Q}$. 

\end{proof}

This completes the proof of Theorem \ref{thm:main}. 

\end{proof}

\section{Open Questions}

We consider the next logical questions  in Shafer's program \cite{Sh2} of comparing and contrasting closed degrees of difficulty in Cantor space and Baire space. First, one can replace Medvedev reducibility with Muchnik:
\begin{qu} Does there exist a closed subset of Baire space not \emph{Muchnik} equivalent to any closed subset of Cantor space?
\end{qu}

Achieving the negative answer in this paper for Medvedev equivalence makes essential use of the fact that one solitary Turing reduction is a continuous map. In the Muchnik case a wholly different technique seems necessary, since a reduction does not require one uniform Turing functional. 

Turning back to Medvedev degrees: our result is that the Medvedev degrees of closed sets in Baire and Cantor space are not all the same. On the other hand, Shafer has shown that the the first-order theories of the Medvedev degrees of closed sets in Baire and Cantor space are recursively isomorphic. So Shafer showed their theories have the same degree of unsolvability and we have shown that they are not exactly the same Medvedev degrees. We'd like to zero in on exactly how similar these structures are. Logical next questions, as mentioned by Shafer:

\begin{qu} Are the Medvedev degrees of closed subsets of Baire space isomorphic (as a lattice), to the Medvedev degrees of closed subset of Cantor space? If not, are they elementarily equivalent?

\end{qu}

Before answering Question 2, we probably need to learn more about the local properties of the Medvedev degrees of closed sets in Cantor and Baire space. In fact, relativization often carries over to the general situation of closed sets the results that hold for effectively closed sets (aka $\Pi_1^0$ classes). Which leads to where we think the next progress will come from. 

In some sense, the next natural question is to decide what sorts of lattice extensions are possible, given any particular finite configurations of Medvedev degrees of $\Pi_1^0$ mass problems. (Essentially the extension of embeddings problem: see Lerman, \cite{Le}.)  Accomplishing this would lead to an affirmative answer to the following question:
\begin{qu} Is the $\forall \exists$-theory in the language of lattices of the Medvedev degrees of 
\newline
$\Pi_1^0$ subsets of Cantor (or Baire) space decidable?
\label{qu:ext}
\end{qu}

In the language of partial orders and for Cantor space the question has been answered in the affirmative by Cole and Kihara \cite{CoKi, Cole}. Their proof built on the proof of the density of the Medvedev degrees of $\Pi_1^0$ subsets of Cantor space (Cenzer/Hinman) \cite{CeHi}.

Sometimes proofs for results about Muchnik degrees give insight into proofs about Medvedev degrees and vice-versa, so another interesting (still open) question is: 

\begin{qu} Is the $\forall \exists$-theory (in the language of partial orders) of the Muchnik degrees of $\Pi_1^0$ subsets of Cantor (or Baire) space decidable?
\label{qu:extMuchnik}
\end{qu}

We hope that the techniques in this paper can be combined with the density construction for the Muchnik degrees of closed sets \cite{BiShSi} to yield to an affirmative answer to Question \ref{qu:extMuchnik}.

 Interestingly, the proof of density for Muchnik degrees of closed subsets of Cantor space made use of some particular subsets of Baire space, which is one reason we believe the techniques in this paper may be useful even when Question \ref{qu:extMuchnik} is restricted to closed subsets of Cantor space.

\bmhead{Acknowledgements}

I thank Peter G. Hinman for invaluable help. Though retired and never having met me, he kindly accepted my request to discuss research ideas. This led to many discussions, and his extensive feedback has contributed much to this paper. I would also like to thank Sean Walsh, whose suggestions about where to look for interesting and solvable open questions led me to the one this paper settles.




\end{document}